\documentclass[12pt]{article}
\usepackage{graphicx}
\usepackage{amssymb,amsmath,amsfonts,amsthm}
\usepackage{url}
\small\normalsize

\newtheorem{thm}{Theorem}[section]
\newtheorem{lem}[thm]{Lemma}        %% lemmas, props, cor, etc
\newcommand{\qitemm}[1]{\noindent\leavevmode\hangindent10mm%
       \noindent\hbox to10mm{#1\hss}\ignorespaces}
\newcommand{\eop}{\qquad\hspace*{\fill} $\blacksquare$}
\newcommand{\openeop}{\qquad\hspace*{\fill} $\square$}

\begin{document}

\title{\bf Toughness and Vertex Degrees}

\author{{\sc  D. Bauer}\\
  {\small \sl Department of Mathematical Sciences} \\[-4pt]
  {\small \sl Stevens Institute of Technology } \\[-4pt]
  {\small \sl Hoboken, NJ 07030, U.S.A.  }
  \and
  {\sc H.J. Broersma\footnote{\,Supported by EPSRC Grant
      EP/F064551/1.\newline
      \mbox{}\hskip20pt Current address: Department of Applied Mathematics,
      University of Twente, P.O. Box 217, 7500 AE Enschede, The
      Netherlands.}}\\
  {\small \sl School of Engineering}\\[-6pt]
  {\small \sl and Computing Sciences}\\[-4pt]
  {\small \sl Durham University}\\[-4pt]
  {\small \sl South Road, Durham DH1 3LE, U.K.}
  \and
  {\sc J. van den Heuvel}\\
  {\small \sl Department of Mathematics}\\[-4pt]
  {\small \sl London School of Economics}\\[-4pt]
  {\small \sl Houghton Street, London WC2A 2AE, U.K.}
  \and
  {\sc N. Kahl}\\
  {\small \sl Department of Mathematics }\\[-6pt]
  {\small \sl and Computer Science}\\[-4pt]
  {\small \sl Seton Hall University}\\[-4pt]
  {\small \sl South Orange, NJ  07079, U.S.A.}
  \and
  {\sc E. Schmeichel}\\
  {\small \sl Department of Mathematics}\\[-4pt]
  {\small \sl San Jos\'e State University}\\[-4pt]
  {\small \sl San Jos\'e, CA  95192, U.S.A.}
}

\date{}        %%  to get a date, remove this
\maketitle

\begin{abstract}
\noindent
We study theorems giving sufficient conditions on the vertex degrees of a
graph~$G$ to guarantee~$G$ is $t$-tough. We first give a best monotone
theorem when $t\ge1$, but then show that for any integer $k\ge1$, a best
monotone theorem for $t=\frac1k\le 1$ requires at least $f(k)\cdot|V(G)|$
nonredundant conditions, where~$f(k)$ grows superpolynomially as
$k\rightarrow\infty$. When $t<1$, we give an additional, simple theorem
for~$G$ to be $t$-tough, in terms of its vertex degrees.
\end{abstract}

\section{Introduction}\label{s1}

We consider only simple graphs without loops or multiple edges. Our
terminology and notation will be standard except as indicated, and a good
reference for any undefined terms or notation is~\cite{West01}. For two
graphs~$G,H$ on disjoint vertex sets, we denote their \emph{union} by
$G\cup H$. The \emph{join $G+H$} of~$G$ and~$H$ is the graph formed from
$G\cup H$ by adding all edges between~$V(G)$ and~$V(H)$.

For a positive integer~$n$, an \emph{$n$-sequence} (or just a
\emph{sequence}) is an integer sequence $\pi=(d_1,d_2,\dots,d_n)$, with
$0\le d_j\le n-1$ for all~$j$. In contrast to~\cite{West01}, we will
usually write the sequence in nondecreasing order (and may make this
explicit by writing $\pi=(d_1\le\dots\le d_n)$). We will employ the
standard abbreviated notation for sequences, e.g., $(4,4,4,4,4,5,5,6)$ will
be denoted $4^5\,5^2\,6^1$. If $\pi=(d_1,\dots,d_n)$ and
$\pi'=(d_1',\dots,d_n')$ are two $n$-sequences, we say \emph{$\pi'$
  majorizes~$\pi$}, denoted $\pi'\ge\pi$, if $d_j'\ge d_j$ for all~$j$.

A \emph{degree sequence} of a graph is any sequence
$\pi=(d_1,d_2,\dots,d_n)$ consisting of the vertex degrees of the graph. A
sequence~$\pi$ is \emph{graphical} if there exists a graph~$G$ having~$\pi$
as one of its degree sequences, in which case we call~$G$ a
\emph{realization} of~$\pi$. If~$P$ is a graph property (e.g., hamiltonian,
$k$-connected, etc.), we call a graphical sequence~$\pi$
\emph{forcibly~$P$} if every realization of~$\pi$ has property~$P$.

\medskip
Historically, the degree sequence of a graph has been used to provide
sufficient conditions for a graph to have certain properties, such as
hamiltonicity or $k$-connectivity. In particular, sufficient conditions
for~$\pi$ to be forcibly hamiltonian were given by several authors,
culminating in the following theorem of Chv\'atal~\cite{Chvatal72}.

\begin{thm}[\cite{Chvatal72}]\label{thm:chvatal}
  \;Let $\pi=(d_1\le\dots\le d_n)$ be a graphical sequence, with $n\ge3$.
  If $d_i\le i<\frac12n$ implies $d_{n-i}\ge n-i$, then~$\pi$ is forcibly
  hamiltonian.
\end{thm}

Unlike its predecessors, Chv\'atal's theorem has the property that if it
does not guarantee that~$\pi$ is forcibly hamiltonian because the condition
fails for some $i<\frac12n$, then~$\pi$ is majorized by
$\pi'=i^i\,(n-i-1)^{n-2i}\,(n-1)^i$, which has a unique nonhamiltonian
realization $K_i+(\overline{K_i}\cup K_{n-2i})$. As we will see below, this
implies that Chv\'atal's theorem is the strongest of an entire class of
theorems giving sufficient degree conditions for~$\pi$ to be forcibly
hamiltonian.

Sufficient conditions for~$\pi$ to be forcibly $k$-connected were given by
several authors, culminating in the following theorem of
Bondy~\cite{Bondy69} (though the form in which we present it is due to
Boesch~\cite{Boesch74}).

\begin{thm}[\cite{Boesch74,Bondy69}]\label{BB}
  \;Let $\pi=(d_1\le\dots\le d_n)$ be a graphical sequence with $n\ge2$,
  and let $1\le k\le n-1$. If $d_i\le i+k-2$ implies $d_{n-k+1}\ge n-i$,
  for $1\le i\le\frac12(n-k+1)$, then~$\pi$ is forcibly $k$-connected.
\end{thm}

Boesch~\cite{Boesch74} also observed that Theorem~\ref{BB} is the strongest
theorem giving sufficient degree conditions for~$\pi$ to be forcibly
$k$-connected, in exactly the same sense as Theorem~\ref{thm:chvatal}.

\medskip
Let~$\omega(G)$ denote the number of components of a graph~$G$. For
$t\ge0$, we call~$G$ \emph{$t$-tough} if $t\cdot\omega(G-X)\le|X|$, for
every $X\subseteq V(G)$ with $\omega(G-X)>1$. The \emph{toughness} of~$G$,
denoted~$\tau(G)$, is the maximum $t\ge0$ for which $G$ is $t$-tough
(taking $\tau(K_n)=n-1$, for all $n\ge1$). So if~$G$ is not complete, then
$\tau(G)=\min\Bigl\{\dfrac{|X|}{\omega(G-X)}\Bigm|\text{$X\subseteq V(G)$
  is a cutset of $G$}\Bigr\}$.

In this paper we consider forcibly $t$-tough theorems, for any $t\ge0$.
When trying to formulate and prove this type of theorem, we encountered
very different behavior in the number of conditions required for a best
possible theorem for the cases $t\ge1$ and $t<1$. In order to describe this
behavior precisely, we need to say what we mean by a `condition' and by a
`best possible theorem'.

\medskip
First note that the conditions in Theorems~\ref{thm:chvatal} can be written
in the form:
\[\text{$d_i\ge i+1$ \ or \ $d_{n-i}\ge n-i$, \ for
  $i=1,\ldots,\bigl\lfloor\tfrac12(n-1)\bigr\rfloor$,}\]
and the conditions in Theorem~\ref{BB} can be written in a similar way. We
will use the term `Chv\'atal-type conditions' for such conditions.
Formally, a \emph{Chv\'atal-type condition} for $n$-sequences
$(d_1\le d_2\le\cdots\le d_n)$ is a condition of the form
\[d_{i_1}\ge k_{i_1}\;\vee\;d_{i_2}\ge k_{i_2}\;\vee\;\ldots\;\vee\;
d_{i_r}\ge k_{i_r},\]
where all~$i_j$ and~$k_{i_j}$ are integers, with
$1\le i_1<i_2<\dots<i_r\le n$ and
$1\le k_{i_1}\le k_{i_2}\le\cdots\le k_{i_r}\le n$.

\medskip
A graph property~$P$ is called \emph{increasing} if whenever a graph~$G$
has~$P$, so does every edge-augmented supergraph of~$G$. In particular,
``hamiltonian'', ``$k$-connected'' and ``\mbox{$t$-tough}'' are all
increasing graph properties. In this paper, the term ``graph property''
will always mean an increasing graph property.

Given a graph property~$P$, consider a theorem~$T$ which declares certain
degree sequences to be forcibly~$P$, rendering no decision on the remaining
degree sequences. We call such a theorem~$T$ a \emph{forcibly $P$-theorem}
(or just a \emph{$P$-theorem}, for brevity). Thus Theorem~\ref{thm:chvatal}
would be a forcibly hamiltonian theorem. We call a $P$-theorem~$T$
\emph{monotone} if, for any two degree sequences $\pi,\pi'$, whenever~$T$
declares~$\pi$ forcibly~$P$ and $\pi'\ge\pi$, then~$T$ declares~$\pi'$
forcibly~$P$. We call a $P$-theorem~$T$ \emph{optimal} if whenever~$T$ does
not declare a degree sequence~$\pi$ forcibly~$P$, then~$\pi$ is not
forcibly~$P$; $T$ is \emph{weakly optimal} if for any sequence~$\pi$ (not
necessarily graphical) which~$T$ does not declare forcibly~$P$, $\pi$ is
majorized by a degree sequence which is not forcibly~$P$.

A $P$-theorem which is both monotone and weakly optimal is a best monotone
\mbox{$P$-theorem}, in the following sense.

\begin{thm}\label{wo}
  \quad Let $T$,~$T_0$ be monotone $P$-theorems, with~$T_0$ weakly optimal.
  If~$T$ declares a degree sequence~$\pi$ to be forcibly~$P$, then so
  does~$T_0$.
\end{thm}

\textbf{Proof of Theorem \ref{wo}:} Suppose to the contrary that there
exists a degree sequence~$\pi$ so that~$T$ declares~$\pi$ forcibly~$P$,
but~$T_0$ does not. Since~$T_0$ is weakly optimal, there exists a degree
sequence $\pi'\ge\pi$ which is not forcibly~$P$. This means that also~$T$
will not declare~$\pi'$ forcibly~$P$. But if~$T$ declares~$\pi$
forcibly~$P$, $\pi'\ge\pi$, and~$T$ does not declare~$\pi'$ forcibly~$P$,
then~$T$ is not monotone, a contradiction.\eop

\bigskip
If~$T_0$ is Chv\'atal's hamiltonian theorem (Theorem~\ref{thm:chvatal}),
then~$T_0$ is clearly monotone, and we noted above that~$T_0$ is weakly
optimal. So by Theorem~\ref{wo}, Chv\'atal's theorem is a best monotone
hamiltonian theorem.

\medskip
Our goal in this paper is to consider forcibly $t$-tough theorems, for any
$t\ge0$. In Section~\ref{s2} we first give a best monotone $t$-tough
theorem for $n$-sequences, requiring at most
$\bigl\lfloor\frac12n\bigr\rfloor$ Chv\'atal-type conditions, for any
$t\ge1$. In contrast to this, in Sections~\ref{sub1.1} and~\ref{s2a} we
show that for any integer $k\ge1$, a best monotone $1/k$-tough theorem
contains at least $f(k)\cdot n$ nonredundant Chv\'atal-type conditions,
where~$f(k)$ grows superpolynomially as $k\rightarrow\infty$. A similar
superpolynomial growth in the complexity of the best monotone
$k$-edge-connected theorem in terms of~$k$ was previously noted by
Kriesell~\cite{Kp}.

This superpolynomial complexity of a best monotone $1/k$-tough theorem
suggests the desirability of finding more reasonable $t$-tough theorems,
when $t<1$. In Section~\ref{s3} we give one such theorem. This theorem is a
monotone, though not best monotone, $t$-tough theorem which is valid for
any $t\le1$.

\section{A Best Monotone $t$-Tough Theorem for $t\ge1$}\label{s2}

We first give a best monotone $t$-tough theorem for $t\ge1$.

\begin{thm}\label{tge1}
  \quad Let $t\ge1$, $n\ge\lceil t\rceil+2$, and let $\pi=(d_1\le\dots\le
  d_n)$ be a graphical sequence. If{

    \qitemm{($*t$)}\hspace*{\fill}%
    $d_{\lfloor i/t\rfloor}\ge i+1$ \ or \
    $d_{n-i}\ge n-\lfloor i/t\rfloor$, \ for
    $t\le i<\dfrac{tn}{(t+1)}$,\hspace*{\fill}

  }then~$\pi$ is forcibly $t$-tough.
\end{thm}

Clearly, property~($*t$) in Theorem~\ref{tge1} is monotone. Furthermore,
if~$\pi$ does not satisfy~($*t$) for some~$i$ with $t\le i<tn/(t+1)$,
then~$\pi$ is majorized by
$\pi'=i^{\lfloor i/t\rfloor}\linebreak[1]\,
(n-\lfloor i/t\rfloor-1)^{n-i-\lfloor i/t\rfloor}\,(n-1)^i$, which has the
non-$t$-tough realization
$K_i+\bigl(\overline{K_{\lfloor i/t\rfloor}}\cup
K_{n-i-\lfloor i/t\rfloor}\bigr)$. Thus~($*t$) in Theorem~\ref{tge1} is
also weakly optimal, and so Theorem~\ref{tge1} is best monotone by
Theorem~\ref{wo}. Finally, note that when $t=1$, ($*t$) reduces to
Chv\'atal's hamiltonian condition in Theorem~\ref{thm:chvatal}.

\bigskip
\textbf{Proof of Theorem \ref{tge1}:} Suppose~$\pi$ satisfies~($*t$) for
some $t\ge1$ and $n\ge\lceil t\rceil+2$, but~$\pi$ has a realization~$G$
which is not $t$-tough. Then there exists a set $X\subseteq V(G)$ that is
maximal with respect to $\omega(G-X)\ge2$ and $\dfrac{|X|}{\omega(G-X)}<t$.
Let $x\doteq |X|$ and $w\doteq\omega(G-X)$, so that
$w\ge\lfloor x/t\rfloor+1$. Also, let $H_1,H_2,\dots,H_w$ denote the
components of $G-X$, with $|H_1|\ge|H_2|\ge\cdots\ge|H_w|$, and let
$h_j\doteq|H_j|$ for $j=1,\ldots,w$. By adding edges (if needed) to~$G$, we
may assume $\langle X\rangle$ is complete, and each $\langle H_j\rangle$ is
complete and completely joined to~$X$.

Set $i\doteq x+h_2-1$.

\bigskip
\textbf{Claim 1.}\quad$i\ge t$.

\medskip
\textbf{Proof:} It is enough to show that $x\ge t$. Assume instead that
$x<t$. Define $X'\doteq X\cup\{v\}$, with $v\in H_1$. If $h_1\ge2$, then
\[\frac{|X'|}{\omega(G-X')}\:=\:\frac{x+1}{\omega(G-X)}\:<\:\frac{t+1}2\:
\le\:t,\]
which contradicts the maximality of~$X$. Similarly, if $h_1=1$ and $w\ge3$,
then
\[\frac{|X'|}{\omega(G-X')}\:=\:\frac{x+1}{\omega(G-X)-1}\:<\:\frac{t+1}2\:
\le\:t,\]
also a contradiction. Finally, if $h_1=1$ and $w=2$, then~$G$ is the graph
$K_{n-2}+\overline{K_2}$ with $n-2=x<t$, contradicting
$n\ge\lceil t\rceil+2$.\openeop

\bigskip
\textbf{Claim 2.}\quad$i<\dfrac{tn}{t+1}$

\medskip
\textbf{Proof:} Note that $n=x+h_1+h_2+\cdots+h_w\ge x+2h_2+w-2$. Since
$x<tw$, we obtain
\begin{align*}
  i\:&=\:x+h_2-1\:=\:
  \frac{tx+x+(t+1)(h_2-1)}{t+1}\\[1mm]
  &<\:\frac{t(x+w+(1+1/t)(h_2-1))}{t+1}\:\le\:\frac{t(x+2h_2+w-2)}{t+1}
  \:\le\:\frac{tn}{t+1}.
\end{align*}

\vspace*{-9mm}
\openeop

\medskip
\bigskip
By the claims we have $t\le i<\dfrac{tn}{t+1}$. Next note that
\[d_{\lfloor i/t\rfloor}\:=\:d_{\lfloor(x+h_2-1)/t\rfloor}\:\le\:
  d_{\lfloor x/t\rfloor+h_2-1}\:\le\:d_{w+h_2-2}\:\le\:d_{(h_2+\cdots+h_w)}
  \:=\:x+h_2-1\:=\:i.\]
However, we also have
\begin{align*}
d_{n-i}\:&\le\:d_{n-x}\:=\:x+h_1-1\:=\:n-h_2-(h_3+\cdots+h_w)-1\:\le\:
n-(w+h_2-1)\\[1mm]
&<\:n-\Bigl(\frac{x}{t}+h_2-1\Bigr)\:\le\:
n-\frac{x+h_2-1}{t}\:=\:n-i/t\:\le\:n-\lfloor i/t\rfloor,
\end{align*}
contradicting~($*t$).\eop

\section{The Number of Chv\'atal-Type Conditions in\\ Best Monotone
  Theorems}\label{sub1.1}

In this section we provide a theory that allows us to lower bound the
number of degree sequence conditions required in a best monotone
$P$-theorem.

Recall that a \emph{Chv\'atal-type condition} for $n$-sequences
$(d_1\le d_2\le\cdots\le d_n)$ is a condition of the form
\[d_{i_1}\ge k_{i_1}\;\vee\;d_{i_2}\ge k_{i_2}\;\vee\;\ldots\;\vee\;
d_{i_r}\ge k_{i_r},\]
where all~$i_j$ and~$k_{i_j}$ are integers, with
$1\le i_1<i_2<\dots<i_r\le n$ and
$1\le k_{i_1}\le k_{i_2}\le \cdots\le k_{i_r}\le n$. Given an $n$-sequence
$\pi=(k_1\le k_2\le\cdots\le k_n)$, let $C(\pi)$ denote the Chv\'atal-type
condition:
\[d_1\ge k_1+1\;\vee\;d_2\ge k_2+1\;\vee\:\ldots\;\vee\;d_n\ge k_n+1.\]
Intuitively, $C(\pi)$ is the weakest condition that `blocks'~$\pi$. For
instance, if $\pi=2^23^35$, then~$C(\pi)$ is
\begin{equation}\label{eq11.1}
  d_1\ge3\;\vee\;d_2\ge3\;\vee\;d_3\ge4\;\vee\;d_4\ge4\;\vee\;d_5\ge4\;
\vee\;d_6\ge6.
\end{equation}
Since $n$-sequences are assumed to be nondecreasing, $d_1\ge3$ implies
$d_2\ge3$, etc. Also, we cannot have $d_i\ge n$, so the condition $d_6\ge6$
is redundant. Hence~\eqref{eq11.1} can be simplified to the equivalent
Chv\'atal-type condition
\begin{equation}\label{eq11.2}
  d_2\ge3\;\vee\;d_5\ge4,
\end{equation}
and we use $\eqref{eq11.1}\cong\eqref{eq11.2}$ to denote this equivalence.

Conversely, given a Chv\'atal-type condition~$c$, let~$\Pi(c)$ denote the
minimal $n$-sequence that majorizes all sequences which violate~$c$
($\Pi(c)$ may not be graphical). So if~$c$ is the condition
in~\eqref{eq11.2} and $n=6$, then $\Pi(c)$ is $2^23^35$. Of course,
$\Pi(c)$ itself violates~$c$. Note that~$C$ and~$\Pi$ are inverses: For any
Chv\'atal-type condition~$c$ we have $C(\Pi(c))\cong c$, and for any
$n$-sequence~$\pi$ we have $\Pi(C(\pi))=\pi$.

\medskip
Given a graph property~$P$, we call a Chv\'atal-type degree condition~$c$
\emph{$P$-weakly-optimal} if any sequence~$\pi$ (not necessarily graphical)
which does not satisfy~$c$ is majorized by a degree sequence which is not
forcibly~$P$. In particular, each of the
$\bigl\lfloor\frac12(n-1)\bigr\rfloor$ conditions in Chv\'atal's
hamiltonian theorem is weakly optimal.

\medskip
Next consider the poset whose elements are the graphical sequences of
length~$n$, with the majorization relation $\pi\le\pi'$ as the partial
order relation. We call this poset the \emph{$n$-degree-poset}. Posets of
integer sequences with a different order relation were previously used by
Aigner \& Triesch~\cite{AT} in their work on graphical sequences.

Given a graph property~$P$, consider the set of $n$-vertex graphs without
property~$P$ which are edge-maximal in this regard. The degree sequences of
these edge-maximal, \mbox{non-$P$} graphs induce a subposet of the
$n$-degree-poset, called the \emph{$P$-subposet}. We refer to the maximal
elements of this $P$-subposet as \emph{sinks}, and denote their number by
$s(n,P)$.

\medskip
We first prove the following lemma.

\begin{lem}\label{lem1}
  \quad Let~$P$ be a graph property. If a sink~$\pi$ of the $P$-subposet
  violates a $P$-weakly-optimal Chv\'atal-type condition~$c$, then
  $c\cong C(\pi)$.
\end{lem}

\textbf{Proof:} Since~$\pi$ violates~$c$, $\pi\le\Pi(c)$. Since~$\Pi(c)$
violates~$c$, and~$c$ is $P$-weakly-optimal, there is a sequence
$\pi'\ge\Pi(c)$ such that~$\pi'$ has a non-$P$ realization. But
$\pi'\le\pi''$ for some sink~$\pi''$, giving $\pi\le\Pi(c)\le\pi'\le\pi''$.
Since distinct sinks are incomparable, $\pi=\pi''$. This implies
$\Pi(c)=\pi$, and thus $c\cong C(\Pi(c))\cong C(\pi)$.\eop

\begin{thm}\label{poset}
  \quad Let~$P$ be a graph property. Then any \mbox{$P$-theorem} for
  $n$-sequences whose hypothesis consists solely of $P$-weakly-optimal
  Chv\'atal-type conditions must contain at least $s(n,P)$ such conditions.
\end{thm}

\textbf{Proof:} Consider a $P$-theorem whose hypothesis consists solely of
$P$-weakly-optimal Chv\'atal-type conditions. By Lemma~\ref{lem1}, a
sink~$\pi$ satisfies every Chv\'atal-type condition besides~$C(\pi)$. So
the theorem must include all the Chv\'atal-type conditions~$C(\pi)$,
as~$\pi$ ranges over the $s(n,P)$ sinks.\eop

\bigskip
On the other hand, it is easy to see that if we take the collection of
Chv\'atal-type conditions~$C(\pi)$ for all sinks~$\pi$ in the $P$-subposet,
then this gives a best monotone $P$-theorem.

\medskip
We do not have a comparable result for $P$-theorems if we do not require
the conditions to be $P$-weakly-optimal, let alone if we consider
conditions that are not of Chv\'atal-type. On the other hand, all results
we have discussed so far, and most of the forcibly $P$-theorems we know in
the literature, involve only $P$-weakly-optimal Chv\'atal-type degree
conditions.

\section{Best Monotone $t$-Tough Theorems for $t\le1$}\label{s2a}

Using the terminology from Section~\ref{sub1.1}, it follows that
Theorem~\ref{tge1} gives, for $t\ge1$, a best monotone $t$-tough theorem
using a linear number (in~$n$) of weakly optimal Chv\'atal-type conditions.
On the other hand, we now show that for any integer $k\ge1$, a best
monotone \mbox{$1/k$-tough} theorem for $n$-sequences requires at least
$f(k)\cdot n$ weakly optimal Chv\'atal-type conditions, where~$f(k)$ grows
superpolynomially as $k\rightarrow\infty$. In view of Theorem~\ref{poset},
to prove this assertion it suffices to prove the following lemma.

\begin{lem}\label{superpoly}\quad
  Let $k\ge2$ be an integer, and let $n=m(k+1)$ for some integer $m\ge9$.
  Then the number of ($1/k$-tough)-subposet sinks in the
  $n$-degree-subposet is at least $\dfrac{p(k-1)}{5(k+1)}n$, where~$p$
  denotes the integer partition function.
\end{lem}

Recall that the integer partition function~$p(r)$ counts the number of ways
a positive integer~$r$ can be written as a sum of positive integers. Since
$p(r)\sim\dfrac1{4r\sqrt3}e^{\pi\sqrt{2r/3}}$ as $r\rightarrow\infty$
\cite{HR18}, $f(k)=\dfrac{p(k-1)}{5(k+1)}$ grows superpolynomially as
$k\rightarrow\infty$.

\bigskip
\textbf{Proof of Lemma \ref{superpoly}:} Consider the
collection~$\mathcal{C}$ of all connected graphs on~$n$ vertices which are
edge-maximally not-($1/k$-tough). Each $G\in\mathcal{C}$ has the form
$G=K_j+(K_{c_1}\cup\dots\cup K_{c_{kj+1}})$, where $j<n/(k+1)=m$, so that
$1\le j\le m-1$, and $c_1+\dots+c_{kj+1}$ is a partition of $n-j$. Assuming
$c_1\le\dots\le c_{kj+1}$, the degree sequence of~$G$ becomes
$\pi\doteq(c_1+j-1)^{c_1}\,\dots\,(c_{kj+1}+j-1)^{c_{kj+1}}\,(n-1)^j$. Note
that~$\pi$ cannot be majorized by the degrees of any disconnected graph
on~$n$ vertices, since a disconnected graph has no vertex of degree $n-1$.
By a \emph{complete degree} of a degree sequence we mean an entry in the
sequence equal to $n-1$.

Partition the degree sequences of the graphs in~$\mathcal{C}$ into $m-1$
groups, where the sequences in the $j^\text{th}$ group, $1\le j\le m-1$,
are precisely those containing~$j$ complete degrees. We establish two basic
properties of the $j^\text{th}$ group.

\bigskip
\textbf{Claim 1.}\quad\textit{There are exactly
  $p_{kj+1}\bigl((k+1)(m-j)-1\bigr)$ sequences in the $j^\text{th}$ group.}

\medskip
Here $p_{\ell}(r)$ denotes the number of partitions of integer~$r$ into at
most~$\ell$ parts, or equivalently the number of partitions of~$r$ with
largest part at most~$\ell$.

\medskip
\textbf{Proof of Claim 1:} Each sequence in the $j^\text{th}$ group
corresponds uniquely to a set of $kj+1$ component sizes which sum to $n-j$.
If we subtract~1 from each of those component sizes, we obtain a
corresponding collection of $kj+1$ integers (some possibly~0) which sum to
$n-j-(kj+1)=(k+1)(m-j)-1$, and which therefore form a partition of
$(k+1)(m-j)-1$ into at most $kj+1$ parts.\openeop

\bigskip
\textbf{Claim 2.}\quad \textit{No sequence in the $j^\text{th}$ group
  majorizes another sequence in the $j^\text{th}$ group.}

\medskip
\textbf{Proof:} Suppose the sequences
$\pi\doteq(c_1+j-1)^{c_1}\,\dots\,(c_{kj+1}+j-1)^{c_{kj+1}}\,(n-1)^j$ and
$\pi'\doteq(c'_1+j-1)^{c'_1}\,\dots\,(c'_{kj+1}+j-1)^{c'_{kj+1}}\,(n-1)^j$
are in the $j^\text{th}$ group, with $\pi\ge\pi'$. Deleting the~$j$
complete degrees from each sequence gives sequences
$\sigma\doteq(c_1-1)^{c_1}\,\dots\,(c_{kj+1}-1)^{c_{kj+1}}$ and
$\sigma'\doteq(c'_1-1)^{c'_1}\,\dots\,(c'_{kj+1}-1)^{c'_{kj+1}}$, with
$\sigma\ge\sigma'$.

Let~$m$ be the smallest index with $c_m\ne c'_m$; since $\sigma\ge\sigma'$,
we have $c_m>c'_m$. In particular, $c_1+\dots+c_m>c'_1+\dots+c'_m$. But
$c_1+\dots+c_{kj+1}=c'_1+\dots+c'_{kj+1}=n-j$, and so there exists a
smallest index $\ell>m$ with $c_1+\dots+c_\ell\le c'_1+\dots+c'_\ell$. In
particular, $c_\ell<c'_\ell$. Since
$c'_1+\dots+c'_{\ell-1}<c_1+\dots+c_{\ell-1}<c_1+\dots+c_\ell\le
c_1+\dots+c'_\ell$, we have
$d_{c_1+\dots+c_\ell}=c_\ell-1<c'_\ell-1=d'_{c_1+\dots+c_\ell}$, and thus
$\sigma\ngeq\sigma'$, a contradiction.\openeop

\bigskip
Since $K_j+(K_{c_1}\cup\dots\cup K_{c_{kj+1}})$ has~$n$ vertices,
$K_{c_{kj+1}}$ has at most $n-j-kj$ vertices. This means the largest
possible noncomplete degree in a sequence in the $j^\text{th}$ group is
$j+(n-j-kj-1)=n-kj-1$. Using this observation we can prove the following.

\bigskip
\textbf{Claim 3.}\quad\textit{If a sequence
  $\pi=\cdots\,d^{d-j+1}\,(n-1)^j$ in the $j^\text{th}$ group has largest
  noncomplete degree $d\ge n-k(j+1)$, then~$\pi$ is not majorized by any
  sequence in the $i^\text{th}$ group, for $i\ge j+1$.}

\medskip
In particular, such a~$\pi$ is a sink, since~$\pi$ is certainly not
majorized by another sequence in the $j^\text{th}$ group by Claim~2, nor by
a sequence in groups $1,2,\dots,j-1$, since any such sequence has fewer
than~$j$ complete degrees.

\medskip
\textbf{Proof of Claim 3:} If $d\ge n-k(j+1)$, then the $d+1$ largest
degrees $d^{d-j+1}\,(n-1)^j$ in~$\pi$ could be majorized only by complete
degrees in a sequence in group $i\ge j+1$, since the largest noncomplete
degree in any sequence in group~$i$ is at most $n-ki-1<n-k(j+1)$. There are
only $i\le m-1$ complete degrees in a sequence in group~$i$. On the other
hand, since $j+1\le i<m$, we have $d+1\ge n-k(j+1)+1>m(k+1)-km+1=m+1>m-1$,
a contradiction.\openeop

\bigskip
So by Claim 3, the sequences~$\pi$ in the $j^\text{th}$ group which could
possibly be nonsinks (i.e., majorized by a sequence in group~$i$, for some
$i\ge j+1$), must have largest noncomplete degree at most $n-k(j+1)-1$. So
in a graph $G\in\mathcal{C}$, $G=K_j+(K_{c_1}\cup\dots\cup K_{c_{kj+1}})$,
which realizes a nonsink~$\pi$, each of the~$K_c$'s must have order at most
$(n-k(j+1)-1)-j+1=(k+1)(m-j)-k$. Subtracting~1 from the order of each of
these components gives a sequence of $kj+1$ integers (some possibly~0)
which sum to $(n-j)-(kj+1)=(k+1)(m-j)-1$, and which have largest part at
most $(k+1)(m-j)-k-1=(k+1)(m-j-1)$. Thus there are exactly
$p_{(k+1)(m-j-1)}\bigl((k+1)(m-j)-1)\bigr)$ such sequences, and so there
are at most this many nonsinks in the $j^\text{th}$ group. Setting
$N(j)\doteq (k+1)(m-j)-1$, so that $(k+1)(m-j-1)=N(j)-k$, this becomes at
most $p_{N(j)-k}\bigl(N(j)\bigr)$ nonsinks in the $j^\text{th}$ group of
sequences.

But by Claim~1, there are exactly $p_{kj+1}\bigl(N(j)\bigr)$ sequences in
group~$j$, and so the number of sinks in the $j^\text{th}$ group is at
least $p_{kj+1}\bigl(N(j)\bigr)-p_{N(j)-k}\bigl(N(j)\bigr)$.

Note that $p_{kj+1}(N(j))$ reduces to $p(N(j))$ if $kj+1\ge N(j)$. However,
$kj+1\ge N(j)$ is equivalent to $j\ge\dfrac{(k+1)m-2}{2k+1}$. Since
$k\ge2$, the inequality $j\ge\dfrac{(k+1)m-2}{2k+1}$ holds if
$j\ge\frac35m$. Thus $p_{kj+1}(N(j))=p(N(j))$ holds for $j\ge\frac35m$.

On the other hand, for $j\le m-2$ we can show the following.

\bigskip
\textbf{Claim 4.}\quad\textit{If $j\le m-2$, then
  \[p\bigl(N(j)\bigr)-p_{N(j)-k}\bigl(N(j)\bigr)\:=\:
  1+p(1)+\dots+p(k-1)\:\ge\:p(k-1).\]}%

\textbf{Proof:} Note that if $j\le m-2$, then $k<\frac12N(j)$. The left
side of the equality in the claim counts partitions of~$N(j)$ with largest
part at least $N(j)-(k-1)$. The right side counts the same according to the
exact order $N(j)-\ell$, $0\le\ell\le k-1$, of the largest part in the
partition, using that the largest part is unique since
$N(j)-\ell\ge N(j)-(k-1)>\frac12N(j)$.\openeop

\bigskip
Completing the proof of Lemma~\ref{superpoly}, we find that the number of
sinks in the \mbox{($1/k$-tough)}-subposet of the $n$-degree-poset is at
least
\begin{align*}
  \sum_{j=\lceil 3m/5\rceil}^{m-2}\bigl[p_{kj+1}\bigl(N(j)\bigr)-
  p_{N(j)-k}\bigl(N(j)\bigr)\bigr]\:&=\:
  \sum_{j=\lceil3m/5\rceil}^{m-2}\bigl[p\bigl(N(j)\bigr)-
  p_{N(j)-k}\bigl(N(j)\bigr)\bigr]\\
  &\hskip-15mm\ge\:\sum_{j=\lceil3m/5\rceil}^{m-2}p(k-1)\:\ge\:
  \bigl(\tfrac25m-\tfrac95\bigr)p(k-1)\\
  &\hskip-30mm=\:\Bigl(\frac{2n}{5(k+1)}-\frac95\Bigr)p(k-1)\:\ge\:
  \frac{n}{5(k+1)}p(k-1),
\end{align*}
as asserted, since $n=m(k+1)\ge9(k+1)$ implies
$\dfrac{2n}{5(k+1)}-\dfrac95\ge\dfrac{n}{5(k+1)}$.\eop

\bigskip
Combining Lemma~\ref{superpoly} with Theorem~\ref{poset} gives the promised
superpolynomial growth in the number of weakly optimal Chv\'atal-type
conditions for $1/k$-toughness.

\begin{thm}\label{superthm}\quad
  Let $k\ge2$ be an integer, and let $n=m(k+1)$ for some integer $m\ge9$.
  Then a best monotone $1/k$-tough theorem for $n$-sequences whose degree
  conditions consist solely of weakly optimal Chv\'atal-type conditions
  requires at least $\dfrac{p(k-1)n}{5(k+1)}$ such conditions, where $p(r)$
  is the integer partition function.
\end{thm}

\section{A Simple $t$-Tough Theorem}\label{s3}

The superpolynomial complexity as $k\rightarrow\infty$ of a best monotone
$1/k$-tough theorem suggests the desirability of finding simple $t$-tough
theorems, when $t<1$. We give such a theorem below. It will again be
convenient to assume at first that $t=1/k$, for some integer $k\ge1$. Note
that the conditions in the theorem are still Chv\'atal-type conditions.

\clearpage
\begin{lem}\label{1/kthm}\quad
  Let $k\ge1$ be an integer, $n\ge k+2$, and $\pi=(d_1\le\dots\le d_n)$ a
  graphical sequence. If{

    \qitemm{(i)}$d_i\ge i-k+2$ \ or \ $d_{n-i+k-1}\ge n-i$, \ for
    $k\le i<\frac12(n+k-1)$, and

    \qitemm{(ii)} $d_i\ge i$ \ or \ $d_n\ge n-i$, \ for
    $1\le i\le\frac12n$,

  }then $\pi$ is forcibly $1/k$-tough.
\end{lem}

\textbf{Proof of Lemma~\ref{1/kthm}:} Suppose~$\pi$ has a realization~$G$
which is not $1/k$-tough. By~(ii) and Theorem~\ref{BB}, $G$ is connected.
So we may assume (by adding edges if necessary) that there exists
$X\subseteq V(G)$, with $x\doteq|X|\ge1$, such that $G=K_x+(K_{a_1}\cup
K_{a_2}\cup\dots\cup K_{a_{kx+1}})$, where $1\le a_1\le a_2\le\dots\le
a_{kx+1}$.

Set $i\doteq x+k-2+a_{kx}$.

\bigskip
\textbf{Claim 1.}\quad$k\le i<\frac12(n+k-1)$

\medskip
\textbf{Proof:} The fact that $i\ge k$ follows immediately from the
definition of~$i$. Since $kx-x-k+1=(k-1)(x-1)\ge0$, we have
\begin{equation}\label{eq3}
  kx-1\:\ge\:x+k-2.
\end{equation}
This leads to
\begin{align*}
  n\:=\:x+\sum_{j=1}^{kx-1}a_j+a_{kx}+a_{kx+1}\:&\ge\:x+kx-1+2a_{kx}\\
  &\ge\:2x+k-2+2a_{kx}\:=\:2i-k+2,
\end{align*}
which is equivalent to $i<\frac12(n+k-1)$.\openeop

\bigskip
\textbf{Claim 2.}\quad$d_i\le i-k+1$.

\medskip
\textbf{Proof:} From~\eqref{eq3} we get
\begin{equation}\label{eq4}
  i\:=\:x+k-2+a_{kx}\:\le\:kx-1+a_{kx}\:\le\:\sum_{j=1}^{kx}a_j.
\end{equation}
This gives $d_i\le x+(a_{kx}-1)\:=\:i-k+1$.\openeop

\bigskip
\textbf{Claim 3.}\quad$d_{n-i+k-1}<n-i$.

\medskip
\textbf{Proof:} We have
$n-i+k-1=n-x-a_{kx}+1\le\sum\limits_{j=1}^{kx+1}a_j$. Thus, using the
bound~\eqref{eq4} for~$i$,
\[d_{n-i+k-1}\:\le\:x+a_{kx+1}-1\:<\:n-\sum_{j=1}^{kx}a_j\:\le\:n-i.\]

\vspace*{-9mm}
\openeop

\medskip
\bigskip
Claims~1, 2 and~3 together contradict condition~(i), completing the proof
of the lemma\eop

\bigskip
We can extend Lemma~\ref{1/kthm} to arbitrary $t\le1$ by letting
$k=\lfloor1/t\rfloor$.

\begin{thm}\label{simple1}\quad
  Let $t\le1$, $n\ge\lfloor 1/t\rfloor+2$, and $\pi=(d_1\le\dots\le d_n)$ a
  graphical sequence. If{

    \qitemm{(i)}$d_i\ge i-\lfloor 1/t\rfloor+2$ \ or \
    $d_{n-i+\lfloor 1/t\rfloor-1}\ge n-i$, \ for
    $\lfloor 1/t\rfloor\le i<\frac12\bigl(n+\lfloor 1/t\rfloor-1\bigr)$,
    and

    \qitemm{(ii)} $d_i\ge i$ \ or \ $d_n\ge n-i$, \ for
    $1\le i\le\frac12n$,

  }then $\pi$ is forcibly $t$-tough.
\end{thm}

\textbf{Proof:} Set $k=\lfloor 1/t\rfloor\ge1$. If~$\pi$ satisfies
conditions (i),~(ii) in Theorem~\ref{simple1}, then~$\pi$ satisfies
conditions (i),~(ii) in Lemma~\ref{1/kthm}, and so is forcibly $1/k$-tough.
But $k=\lfloor 1/t\rfloor\le1/t$ means $1/k\ge t$, and so~$\pi$ is forcibly
$t$-tough.\eop

\bigskip
In summary, if $\dfrac1{k+1}<t\le\dfrac1k$ for some integer $k\ge1$, then
Theorem~\ref{simple1} declares~$\pi$ forcibly $t$-tough precisely if
Lemma~\ref{1/kthm} declares~$\pi$ forcibly $1/k$-tough.

\bigskip
\textbf{Acknowledgements.}\\
The authors thank two anonymous referees for comments and suggestions that
greatly improved the structure and clarity of the paper. We also thank
Michael Yatauro for providing the short argument for Claim~2 in the proof
of Theorem~\ref{tge1}.

\end{document}